\documentclass[a4paper,11pt]{article}
\usepackage[utf8]{inputenc}

\usepackage{bm}
\usepackage{amsmath}
\usepackage{amsthm}

\usepackage{mathrsfs}
\usepackage{graphicx}
\usepackage{epsfig, setspace}
\usepackage{pst-plot}
\usepackage{url}
\usepackage{setspace}
\usepackage{graphicx, transparent, color}
\usepackage[ruled]{algorithm2e}
\usepackage{etex}

\DeclareSymbolFont{rsfs}{U}{rsfs}{m}{n}
\DeclareSymbolFontAlphabet{\mathcal}{rsfs}

\theoremstyle{definition}
\theoremstyle{plain}
\usepackage[margin=2.5cm]{geometry}
\usepackage{tikz}
\usetikzlibrary{shapes,calc}
\usepackage{verbatim}
\usepackage{lineno}
\usepackage{listings}
\usepackage{textcomp}

\usepackage{epsfig}
\usepackage{amsmath}
\usepackage{pst-plot}
\usepackage{amssymb}
\usepackage{amsfonts}
\usepackage{lmodern}
\usepackage{euscript}
\usepackage{eufrak}
\usepackage{oldgerm}

\numberwithin{equation}{section}

\usepackage[nottoc]{tocbibind}
\usepackage{color}
\pagecolor{white}
\usepackage{slantsc}
\begin{document}

\title{Construction and comparative study of Euler method with adaptive IQ and IMQ-RBFs}
 
\date{}
\author{
{\sc Samala Rathan\thanks{Email: rathans.math@iipe.ac.in},}
{\sc Deepit Shah}\\[1.5pt]\\
Department of Humanities \& Sciences,\\
Indian Institute of Petroleum and Energy, Visakhapatnam, India-530003
\\[1pt]
}
\maketitle
\begin{abstract}
The fundamental purpose of the present work is to constitute an enhanced Euler method with adaptive inverse-quadratic and inverse-multi-quadratic radial basis function (RBF) interpolation technique to solve initial value problems. These enhanced methods improve the local convergence of numerical solutions by utilizing a free parameter of radial basis functions. Consistency, convergence and stability analysis are provided to support our claims for each method.  Numerical results show that the accuracy and rate of convergence of each proposed method are the same as the original Euler method or improved by making the local truncation error vanish; thus, the adaptive methods are optimal. \end{abstract}
{\small \textbf{Keywords:} Finite Difference Method, Radial Basis Interpolation, Initial Value Problem, Order of accuracy, Rate of Convergence.\\[0.5pt] }
\section{Introduction}
\label{sec:intro}
When numerically solving initial value problems (IVPs) for first-order ordinary differential equations, classical linear finite difference methods, such as Euler's method, are first considered \cite{Numer}. This linear approach yields a fixed order of accuracy as it is designed based on polynomial interpolation; that is,  the accuracy order is calculated regardless of  smoothness of the solution. If we know or can estimate the derivatives of the solution to a specific order, this information could be used to improve the local accuracy \cite{Butcher}. However, the conventional approach of polynomial interpolation lacks the same adaptability and versatility. Thus, this motivates us to study the radial basis function (RBF) interpolation-based modification of classical methods to improve their order and rate of convergence.
\par
In \cite{GJ-1, GJ-2}, the adaptive RBF solvers have been developed for solving ODEs such as Euler's method, the midpoint method, and the Adams methods by replacing the polynomial basis with multi-quadratic and Gaussian RBFs respectively. The critical concept of adaptive RBF solvers is that the radial basis functions are defined with a free (shape) parameter which can be adjusted according to the local conditions of the solution. There is a possibility to determine such free parameter based on the local smoothness of the solution by substitution of the same for the polynomial basis, that is,  derivatives, in order to terminate the leading truncation error terms.  This idea can be seen in the applications of ENO and WENO for solving hyperbolic PDEs \cite{GJ-3, GJ-4} using a non polynomial basis such as radial basis functions. The ENO and WENO finite difference methods are high-order accurate finite difference methods for problems with piecewise smooth solutions that contain discontinuities. The solution must be evaluated at each cell boundary in both methods, and the standard reconstruction is based on Lagrange interpolation \cite{Shu}.  Various approaches proposed in \cite{F, FW, R, S} on finding the optimal values of the free parameter for the RBF approximation to get the accuracy or stability.  However, there is no complete theory on how to choose the `right' shape parameter value. The role is often decided by the essence of the issue under consideration. 

\par
The local truncation error analysis is essential for the numerical solution of IVPs, and our approach is to focus on the local truncation error. The best shape parameter removes the leading truncation error term(s) and meets the consistency requirement. Our adaptive shape parameter is optimal in this regard. However, it would be more appealing if there was another way to estimate the optimum shape parameter automatically, but such an automated algorithm does not exist.  Including this, development of midpoint method and Adams methods with inverse quadratic and inverse multi-quadratic RBF methods for IVPs  is our future interest.
\par
In Section 2, we begin with adaptive inverse quadratic and inverse multi-quadratic RBF methods for IVPs. First, we develop inverse quadratic and inverse multi-quadratic RBF interpolation. In Section 3, the classical Euler's method is modified to inverse quadratic and inverse multi-quadratic RBF variants. The shape parameter is locally optimized to make the leading error term(s) vanish to increase the order of the local truncation error and further increase the global order of accuracy. We derive the adaptive conditions for enhancing the order of accuracy with higher-order derivatives under the assumption that the solution is smooth. We also evaluate the optimal shape parameter value for higher-order accuracy. In Section 4, we show convergence and stability region for each proposed method. We compare the stability regions of each method i.e., the inverse-quadratic and inverse multi-quadratic version to the original Euler method. In Section 5, numerical results show that inverse quadratic and multi-quadratic RBF methods produce the desired order of accuracy. One method produce better accuracy than the other depending on the method and the problem. We limit our discussion to the case where the RBFs in question has only one shape parameter, though multiple shape parameters may be used for further improvement. The initial finite difference method is the limit case of the corresponding modified RBF method, one of the advantages of modified RBF methods. If the shape parameter is removed, the modified RBF method becomes equal to the original method. Thus, it concludes that even if the shape parameter does not fit or deviates from the optimal value, the modified RBF methods exhibit at least the same convergence rate as the original methods. As a consequence, the order of the modified methods can only get better. Finally, in Section 6, we give some concluding remarks.

\section{RBF Interpolation}\label{sec:RBF}
Let us consider RBF interpolation in one-dimension. Given $(N+1)$ distinct data points $(x_0, u_0), . . . ,  (x_N, u_N )$  with $u_k$ the value of the unknown function $u(x)$ at $x = x_k$, where $ x \in \mathbb{R}$. We use the RBFs, $\phi_k (x) = \phi \big(|x - x_k|, \epsilon_k \big)$, where  $\epsilon_k$ is a shape parameter, to find an interpolant based on the given $(N+1)$ data points. The value of $\epsilon_k$ can vary over $x_k$. The interpolant $r (x)$ takes the form of a weighted sum of RBFs
\begin{equation}\label{eq1}
r(x)=\sum_{k=0}^{N} \lambda_k\phi \big(|x - x_k |, \epsilon_k \big)
\end{equation}
where $\lambda_k$ are the unknown expansion parameters to be determined.
Using interpolation restraints $$r(x_k)=u_k, k= 0,1,...,N,$$ the expansion coefficients $\lambda_k$ satisfy the following linear system
\begin{eqnarray*}
\begin{pmatrix}
\phi(|x_0 - x_0|, \epsilon_0) & \phi( |x_0 - x_1|, \epsilon_1) &.... & \phi(|x_0 - x_N |, \epsilon_N )\\
\phi(|x_1 - x_0|, \epsilon_0) & \phi( |x_1 - x_1|, \epsilon_1) &....& \phi(|x_1 - x_N |, \epsilon_N )\\
... & ... &.... &....\\
... & ... &.... &....\\
... & ... &.... &....\\
\phi(|x_N - x_0|, \epsilon_0) &  \phi( |x_N - x_1|, \epsilon_1) &....& \phi(|x_N - x_N |, \epsilon_N )
\end{pmatrix}
\begin{pmatrix}
\lambda_0\\\lambda_1\\.\\.\\.\\\lambda_N
\end{pmatrix}=
\begin{pmatrix}
u_0\\u_1\\.\\.\\.\\u_N\end{pmatrix}.
\end{eqnarray*}
Until and unless specified explicitly, we consider all the shape parameters $\epsilon_k$ are same, i.e., $\epsilon_k=\epsilon$, for all $k$. In this paper, we use inverse multi-quadratic (IMQ) and inverse quadratic  (IQ)-RBFs.

\subsection{Inverse Multi-Quadratic RBF Interpolation}
\label{sec: IRBF}
The inverse multi-quadratic radial basis function $\phi_k(x)$ is $\displaystyle\frac{1}{\sqrt{1+\epsilon_k^2(x-x_k)^2}}$.  Consider the interpolant  for $N=1$ in \eqref{eq1} is given by
\begin{equation*}
r(x)=\displaystyle\lambda_0\phi_0(x)+\lambda_1\phi_1(x).
\end{equation*}
Using the interpolation condition $r(x_k)=u_k, k= 0,1,$ the interpolation matrix becomes a symmetric matrix with all diagonal entries  $1$.  Thus, we have
\begin{equation*}
\displaystyle\begin{pmatrix}
1&\displaystyle\frac{1}{\sqrt{1+\epsilon^2h^2}}\\\displaystyle\frac{1}{\sqrt{1+\epsilon^2h^2}}&1
\end{pmatrix}
\begin{pmatrix}
\displaystyle\lambda_0\\\lambda_1
\end{pmatrix}=\displaystyle\begin{pmatrix}
u_0\\u_1
\end{pmatrix},
\end{equation*}
where $h=x_1-x_0.$  Solving for $\lambda_k, k=0,1$, we get
$$\displaystyle\lambda_0=\frac{1+\epsilon^2h^2}{\epsilon^2h^2}\bigg(u_0-\frac{u_1}{\sqrt{1+\epsilon^2h^2}}\bigg),$$
$$\displaystyle\lambda_1=\frac{1+\epsilon^2h^2}{\epsilon^2h^2}\bigg(u_1-\frac{u_0}{\sqrt{1+\epsilon^2h^2}}\bigg).$$
Differentiating the interpolant $r(x)$  with respect to $x$, we get
\begin{equation*}
\frac{d}{dx}r(x)=-\frac{\lambda_0(x-x_0)\epsilon^2}{{(1+\epsilon^2(x-x_0)^2})^\frac{3}{2}}-\frac{\lambda_1(x-x_1)\epsilon^2}{{(1+\epsilon^2(x-x_1)^2})^\frac{3}{2}},
\end{equation*}
and evaluating at $x=x_0$, we have
\begin{equation}\label{eq2}
\frac{d}{dx} r(x_0)=\frac{\lambda_1\epsilon^2h}{{(1+\epsilon^2h^2})^\frac{3}{2}}.
\end{equation}
Substituting the value of $\lambda_1$ in \eqref{eq2}, we have
\begin{equation}\label{eq3}
\frac{d}{dx} r(x_0)=\frac{u_1\sqrt{1+\epsilon^2h^2}-u_0}{(1+\epsilon^2h^2)h},
\end{equation}
and as $\epsilon \to 0$, the equation \eqref{eq3} become
\begin{equation*}
\lim_{\epsilon \to 0}\frac{d}{dx}r(x_0)=\frac{u_1-u_0}{h}.
\end{equation*}
Thus, the method reduces to forward difference formula of $u^{\prime}_x$ at $x=x_0.$

\subsection{Inverse Quadratic RBF Interpolation}
\label{sec: IQRBF}
The inverse quadratic radial basis function is $\phi_k(x)=\displaystyle\frac{1}{1+\epsilon_k^2(x-x_k)^2}$. Then, for $N=1$ in \eqref{eq1}  and using the interpolation condition on the interpolant
 \begin{equation*}
r(x)=\displaystyle\lambda_0\phi_0(x)+\lambda_1\phi_1(x),
\end{equation*}
the interpolation matrix becomes a symmetric matrix with all diagonal entries 1, which is
\begin{equation*}
\begin{pmatrix}
1&\displaystyle\frac{1}{1+\epsilon^2h^2}\\\displaystyle\frac{1}{1+\epsilon^2h^2}&1
\end{pmatrix}
\begin{pmatrix}
\displaystyle\lambda_0\\\lambda_1
\end{pmatrix}=\displaystyle\begin{pmatrix}
u_0\\u_1
\end{pmatrix}.
\end{equation*}
Solving for $\lambda_k, k=0,1$, we get
$$\displaystyle\lambda_0=\frac{(1+\epsilon^2h^2)}{\epsilon^2h^2(2+\epsilon^2h^2)}\bigg((1+\epsilon^2h^2)u_0-u_1\bigg),$$
$$\displaystyle\lambda_1=\frac{(1+\epsilon^2h^2)}{\epsilon^2h^2(2+\epsilon^2h^2)}\bigg((1+\epsilon^2h^2)u_1-u_0\bigg).$$
Differentiating $r(x)$ with respect to $x$, we get
\begin{equation}\label{eq4}
\frac{d}{dx}r(x)=-2\frac{\lambda_0(x-x_0)\epsilon^2}{{(1+\epsilon^2(x-x_0)^2})^2}-2\frac{\lambda_1(x-x_1)\epsilon^2}{{(1+\epsilon^2(x-x_1)^2})^2}.\end{equation}
At $x=x_0$ in \eqref{eq4}, we obtain
\begin{equation}\label{eq5}
\frac{d}{dx} r(x_0)=2\frac{\lambda_1\epsilon^2h}{{(1+\epsilon^2h^2})^2}.\end{equation}
Substituting value of $\lambda_1$ in \eqref{eq5}, we get
\begin{equation}\label{eqiq}
\frac{d}{dx} r(x_0)=\dfrac{2\bigg(u_1(1+\epsilon^2h^2)-u_0\bigg)}{\bigg(1+\epsilon^2h^2\bigg)\bigg(2+\epsilon^2h^2\bigg)h},
\end{equation}
and letting  $\epsilon \to 0$, we get
\begin{equation*}
\lim_{\epsilon \to 0}\frac{d}{dx}r(x_0)=\frac{u_1-u_0}{h}.
\end{equation*}
Note that, again the method reduces to forward difference formula of $u^{\prime}_x$ at $x=x_0$.

\section{Adaptive  RBF Euler Method for IVPs}
\label{sec: ARBF}
We consider the initial value problem of the form
\begin{equation} \label{Eq31}
\dfrac{du}{dt}=f(t,u), \,\, a \leq t \leq b,
\end{equation}
with initial condition 
\begin{equation}
u(a)=u_0,
\end{equation}
where we assume $u(t) \in C^{\infty}[a,b]$ and $f(t,u)$ is a class of $C^{\infty}$ function. We divide the interval $[a,b]$ in uniform way $t_n=a+nh,$ $n=0,1,...,N,$ where $h=\dfrac{b-a}{N}$ is the grid size. In this section, we derive the adaptive inverse multi-quadratic, inverse quadratic RBF Euler methods and their modifications.

\subsection{Consistency of adaptive IMQ-RBF Euler Method}
\label{sec: CIMQRBF}
Substitute $\epsilon=\epsilon_n$, $u_1=u_{n+1}$, $u_0=u_n$ in \eqref{eq3}, we get
\begin{equation} \label{Eq32}
\frac{d}{dx} r(x_0)=\frac{u_{n+1}\sqrt{1+\epsilon_n^2h^2}-u_{n}}{(1+\epsilon_n^2h^2)h}.
\end{equation}
Now, we approximate $\dfrac{du}{dt}$ in \eqref{Eq31} by the first derivative of inverse multi-quadratic RBF interpolant given in \eqref{Eq32} at each discrete mesh point $t_n$, $n=0,1,...,N,$ we obtain
\begin{equation}\label{IMQE-RBF}
u_{n+1}=\frac{(1+\epsilon_n^2h^2)hf_n+u_n}{\sqrt{1+\epsilon_n^2h^2}},
\end{equation}
where $f_n=f(t_n, u_n)$. This method differs from the original Euler method. Local truncation error is defined as the error of the method over a time step, normalized by $h.$ Then the local truncation error, $\tau_n$, of the RBF Euler’s method is
$$\displaystyle \tau_n=\displaystyle \dfrac{u_{n+1}-\dfrac{(1+\epsilon_n^2h^2)hf_n+u_n}{\sqrt{1+\epsilon_n^2h^2}}}{h}.$$
Using Taylor series expansion about a point $t=t_n$ yields
\begin{equation}\label{trunc}
\tau_n=h\bigg(\frac{u_n''}{2}+\frac{u_n}{2}\epsilon_n^2\bigg)+h^2\bigg(\frac{u_n'''}{6}+\frac{u_n'}{2}\epsilon_n^2\bigg)+h^3\bigg(\frac{u_n^4}{24}-\frac{3u_n}{8}\epsilon_n^4\bigg)+h^4\bigg(\frac{u_n^5}{120}-\frac{3u_n'}{8}\epsilon_n^4\bigg)+O(h^5).
\end{equation}
Notice that the modified Euler’s method with the IMQ-RBF still yields the first order accuracy as the leading error term is of $O(h)$. However, the coefficient of the first order term is not uniquely determined as it contains the free parameter $\epsilon_n.$
\par
$O(h^2)$: If we allow the leading error term to be zero, then we can eliminate the first term in the truncation error so that we arrive at  second order of convergence. Thus,
$\displaystyle \dfrac{u_n''}{2}+\frac{u_n}{2}\epsilon_n^2=0,$ yields
\begin{equation}\label{second}
\epsilon_n^2=-\frac{u_n''}{u_n}.
\end{equation}
Note that $\epsilon_n^2$ indexed by $n$, therefore the value $\epsilon_n^2$ varies with $n$ as the solution does.
\par 
$O(h^3)$: If we allow the two leading error terms to be zero, then we could arrive at  third order of convergence. Thus
$$\displaystyle \bigg(\frac{u_n''}{2}+\frac{u_n}{2}\epsilon_n^2\bigg)h+\bigg(\frac{u_n'''}{6}+\frac{u_n'}{2}\epsilon_n^2\bigg)h^2=0,$$ yields
\begin{equation}
\epsilon_n^2=-\displaystyle \dfrac{3u_n''
+hu_n'''}{3(u_n+hu_n')}.
\end{equation}
As step size $h \to 0$, we get
\begin{equation}
\lim_{h \to 0}\epsilon_n^2=\lim_{h \to 0}-\displaystyle \dfrac{3u_n''
+hu_n'''}{3(u_n+hu_n')}
\\
=-\frac{u_n''}{u_n},
\end{equation}
which is identical to the value of $\epsilon_n^2$ in \eqref{second}.
\par 
$O(h^4)$: Similarly to achieve fourth of convergence, we annihilate the first three terms in the truncation error \eqref{trunc}. Thus
$$\displaystyle \bigg(\frac{u_n''}{2}+\frac{u_n}{2}\epsilon_n^2\bigg)h+\bigg(\frac{u_n'''}{6}+\frac{u_n'}{2}\epsilon_n^2\bigg)h^2+\bigg(\frac{u_n^4}{24}-\frac{3u_n}{8}\epsilon_n^4\bigg)h^3=0,$$
produce
\begin{equation}\label{two}
\epsilon_n^2=\frac{6(u_n+hu_n') \pm \sqrt{36(u_n+hu_n')^2+9h^2u_n(12u_n''+4hu_n'''+h^2u_n^4)} }{9h^2u_n}.
\end{equation}
Now, this equation  \eqref{two} provide two values for $\epsilon_n^2$. These are
\begin{equation}
\epsilon_n^{2+}=\frac{6(u_n+hu_n') - \sqrt{36(u_n+hu_n')^2+9h^2u_n(12u_n''+4hu_n'''+h^2u_n^4)} }{9h^2u_n},
\end{equation}
\begin{equation}
\epsilon_n^{2-}=\frac{6(u_n+hu_n') + \sqrt{36(u_n+hu_n')^2+9h^2u_n(12u_n''+4hu_n'''+h^2u_n^4)} }{9h^2u_n}.
\end{equation}
Now, we derive a condition on consistency for $\epsilon_n^{2+}$ and $\epsilon_n^{2-}$ based on the value of $u_n$ is non-negative and non-positive.
Let consider the first case as $ u_n>0, u_n \not = 0.$ Upon checking for consistency,
\begin{equation}
\lim_{h \to 0}\epsilon_n^{2+}=\lim_{h \to 0} \frac{6(u_n+hu_n') - \sqrt{36(u_n+hu_n')^2+9h^2u_n(12u_n''+4hu_n'''+h^2u_n^4)} }{9h^2u_n}
=-\frac{u_n''}{u_n}
\end{equation}
and 
\begin{equation}
\lim_{h \to 0}h^2\epsilon_n^{2-}=\lim_{h \to 0} \frac{6(u_n+hu_n') + \sqrt{36(u_n+hu_n')^2+9h^2u_n(12u_n''+4hu_n'''+h^2u_n^4)} }{9h^2u_n}
=\frac{4}{3},
\end{equation}
Thus   $\tau_n \to 0$ as $h \to 0$ for $\epsilon_n^{2+}$ and  $\tau_n \not \to 0$ as $h \to 0$ for $\epsilon_n^{2-}$. This concludes that the consistency holds for $\epsilon_n^{2+}$ which is similar to the condition \eqref{second} but not for $\epsilon_n^{2-}$.
Similarly, for the second case  $ u_n<0, u_n \not = 0$, we arrive at $\epsilon_n^{2+}$ is not consistent  and $\epsilon_n^{2-}$  is consistent.  We can repeat this procedure for achieving fifth and higher-order accuracy but the algebraic procedure become complicated to determine the optimal shape parameter $\epsilon_n^2$. The conjecture stated in   \cite{GJ-1} that there exists at least one shape parameter for the consistency to any order of accuracy and it is main reason why RBF approximations yields spectral accuracy. We do not provide much details here, interested readers can look at \cite{GJ-1, GJ-2}.  
We have derived modified Euler method \eqref{IMQE-RBF}, however, these have some technical difficulties.
\begin{itemize}
\item The square root is involved in the denominator of \eqref{IMQE-RBF} in the approximation of $u_{n+1}$ and in some extreme cases, this value may become complex but the solution is real.
\item To achieve simplest second order of accuracy, the second derivative is involved in the optimal shape parameter $\epsilon_n^2$ and should be estimated. 
\item The optimal shape parameter $\epsilon_n^2$  is given in fraction and the situation may encounter for denominator in \eqref{second} become zero. 
\end{itemize}
To avoid these situations, we further derive adaptive Euler method based on \eqref{IMQE-RBF} as follows. The equation \eqref{IMQE-RBF} is
\begin{equation}\label{IMQE1-RBF} 
u_{n+1}=\frac{(1+\epsilon_n^2h^2)hf_n+u_n}{\sqrt{1+\epsilon_n^2h^2}}.
\end{equation}
Now to avoid the square roots in the approximation, we use Taylor's expansion on $\dfrac{1}{\sqrt{1+\epsilon_n^2h^2}}$ as 
\begin{equation}\label{TE1}
\dfrac{1}{\sqrt{1+\epsilon_n^2h^2}}\approx 1-\displaystyle \dfrac{\epsilon_n^2h^2}{2}+O(h^4).
\end{equation}
Using \eqref{TE1} in \eqref{IMQE1-RBF} and ignoring higher order terms in expansion, we obtain the modified adaptive Euler method as
\begin{equation}\label{MIMQ-RBF}
u_{n+1}=\bigg(1-\displaystyle \dfrac{\epsilon_n^2h^2}{2}\bigg)\bigg((1+\epsilon_n^2h^2)hf_n+u_n\bigg).
\end{equation}
Local truncation error at point  $t=t_n$ for \eqref{MIMQ-RBF} is
\begin{eqnarray*}
\begin{aligned}
\displaystyle \tau_n= & \displaystyle \dfrac{u_{n+1}-\bigg(1-\dfrac{\epsilon_n^2h^2}{2}\bigg)\bigg((1+\epsilon_n^2h^2)hf_n+u_n\bigg)}{h}\\
= & \bigg(\dfrac{u_n\epsilon_n^2}{2}+\dfrac{u_n''}{2}\bigg)h+\bigg(\dfrac{u_n''\epsilon_n^2}{2}+\dfrac{u_n'''}{6}\bigg)h^2+O(h^3).
\end{aligned}
\end{eqnarray*}
For $O(h^2)$ accuracy, we arrive at condition on the optimal shape parameter as
\begin{equation}
\epsilon_n^2=-\frac{u_n''}{u_n}.
\end{equation}
If we replace the second derivative from the relation $u_n'=f_n$, $u_n''=f_n'$ and estimating the same using backward difference formula $\displaystyle\frac{f_n-f_{n-1}}{h}$, we get
\begin{equation}
\epsilon_n^2=-\dfrac{f_n-f_{n-1}}{hu_n}.
\end{equation}
Therefore, we obtain adaptive modified Euler method is of the form
\begin{equation}
u_{n+1}=\bigg(1+\dfrac{h(f_n-f_{n-1})}{2u_n}\bigg)\bigg(\bigg(1-\dfrac{h^2(f_n-f_{n-1})}{u_n}\bigg)f_n+u_n\bigg).
\end{equation}
Finally, to avoid the situation of denominator in \eqref{second} become zero for the optimal shape parameter, we use following two conditions. The first one is exact limit
\begin{eqnarray}\label{C1}
\text{C1}: \epsilon_n^2=
\begin{cases}
0,\,\,\, & \text{if}\,\,\, u_n=0,\\
-\dfrac{u_n''}{u_n},\,\, & \text{if}\,\,\,\, u_n \neq 0,
\end{cases}
\end{eqnarray}
 and the second one is finite difference limit 
\begin{eqnarray}\label{C2}
\text{C2}: \epsilon_n^2=
\begin{cases}
\text{sgn} \bigg((f_n-f_{n-1}) hu_n \bigg)  L , & \text{if}\,\,\,|u_n| \leq h^p,\\
-\dfrac{f_n-f_{n-1}}{hu_n}, & \text{otherwise,}\,\,\,\, 
\end{cases}
\end{eqnarray}
where $0 \leq L < \infty$ is a nonnegative number,  sgn is the sign function and $p \geq 0.$ Note that the condition C2 is the finite difference approximation of C1. For further details on these conditions, refer \cite{GJ-1, GJ-2}. We propose numerical results based on these two conditions in numerical section.

\subsection{Consistency of adaptive IQ-RBF Euler Method}
\label{sec: CIQRBF}
In following similar lines as consistency analysis for IMQ-RBF Euler method, we get from \eqref{eqiq}, an adaptive Euler method based on inverse quadratic radial basis function as 
\begin{equation}\label{IQ1}
\frac{d}{dx} r(x_0)=\dfrac{2\bigg((1+\epsilon_n^2h^2)u_{n+1}-u_n\bigg)}{h(1+\epsilon_n^2h^2)\big(2+\epsilon_n^2h^2\big)}.
\end{equation}
Thus, the approximation for \eqref{Eq31} from \eqref{IQ1} is
\begin{equation}\label{IQ-RBF}
u_{n+1}=\dfrac{h(1+\epsilon_n^2h^2)(2+\epsilon_n^2h^2)f_n+2u_n}{2(1+\epsilon_n^2h^2)}.
\end{equation}
Then the local truncation error $\tau_n$ of the adaptive IQ-RBF Euler’s method is
\begin{equation}
\displaystyle \tau_n=\displaystyle \dfrac{u_{n+1}-\displaystyle\frac{h(1+\epsilon^2h^2)(2+\epsilon^2h^2)f_n+2u_n}{2(1+\epsilon^2h^2)}}{h}.
\end{equation}
Using Taylor series expansion about point $t=t_n$, we get
\begin{equation}
\tau_n=h\bigg(\dfrac{u_n''}{2}+u_n\epsilon_n^2\bigg)+h^2\bigg(\dfrac{u_n'''}{6}-\dfrac{u_n'}{2}\epsilon_n^2\bigg)+h^3\bigg(\dfrac{u_n^4}{24}-u_n\epsilon_n^4\bigg)+h^4\bigg(\dfrac{u_n^5}{120}-\dfrac{3u_n'}{2}\epsilon_n^4\bigg)+O(h^5).
\end{equation}
Notice that the modified Euler’s method with the IQ-RBF still yields the first order accuracy
as the leading error term is of $O(h)$.
\newline
$O(h^2):$ The optimal value of shape parameter $\epsilon_n^2$ for second order of convergence is given by
\begin{equation}
\epsilon_n^2=-\frac{u_n''}{2u_n}.
\end{equation}
$O(h^3):$ The optimal value of shape parameter $\epsilon_n^2$ for third order of convergence is given by
\begin{equation}
\epsilon_n^2=-\displaystyle \dfrac{\displaystyle \frac{u_n''}{2}+\displaystyle \dfrac{hu_n'''}{6}}{u_n-\displaystyle \dfrac{hu_n''}{2}},
\end{equation}
as $h \to 0$, we get
\begin{equation}
\lim_{h \to 0}\epsilon_n^2=\lim_{h \to 0}-\displaystyle \dfrac{3u_n''+hu_n'''}{3(u_n+hu_n')} =-\frac{u_n''}{2u_n}.
\end{equation}
$O(h^4):$ For fourth order of convergence, we obtain two optimal shape parameter  $\epsilon_n^2$ which are
\begin{equation}
\epsilon_n^{2+}=\displaystyle \dfrac{(u_n-\dfrac{hu_n'}{2})-\sqrt{(u_n-\dfrac{hu_n'}{2})^2+4(h^2u_n)(\dfrac{u_n''}{2}+h\dfrac{u_n'''}{6}+h^2\dfrac{u_n^4}{24}} }{2h^2u_n},
\end{equation}
\begin{equation}
\epsilon_n^{2-}=\displaystyle \dfrac{(u_n-\dfrac{hu_n'}{2})+\sqrt{(u_n-\dfrac{hu_n'}{2})^2+4(h^2u_n)(\dfrac{u_n''}{2}+h\dfrac{u_n'''}{6}+h^2\dfrac{u_n^4}{24}} }{2h^2u_n}.
\end{equation}
Checking the consistency for $ u_n>0, u_n \not = 0$, we get
\begin{equation}
\lim_{h \to 0}\epsilon_n^{2+}=\frac{-u_n''}{2u_n},
\end{equation}
and $\tau_n \to 0$ as $h \to 0$. Thus, consistency holds for $\epsilon_n^{2+}.$
Now, for 
\begin{equation}
\lim_{h \to 0}h^2\epsilon_n^{2-}=1,
\end{equation}
which leads $\tau_n \not \to 0$ as $h \to 0$. Thus, consistency does not holds for $\epsilon_n^{2-}.$ Similiarly, $\epsilon_n^{2+}$ is not consistent for $u_n<0$ and $\epsilon_n^{2-}$ is  consistent for $u_n<0.$
Again, we arrived at  conditions 
\begin{itemize}
\item To achieve at least second order of accuracy, the second derivative is involved in the optimal shape parameter $\epsilon_n^2$ and should be estimated.
\item The optimal shape parameter $\epsilon_n^2$  is given in fraction and the situation may encounter for denominator in \eqref{second} become zero,
\end{itemize}
and to remedy that, we propose modified method for
\begin{equation}\label{IQ3}
u_{n+1}=\frac{h(1+\epsilon_n^2h^2)(2+\epsilon_n^2h^2)f_n+2u_n}{2(1+\epsilon_n^2h^2)}.
\end{equation}
Using Taylor expansion of $\dfrac{1}{1+\epsilon_n^2h^2}\approx 1-\epsilon_n^2h^2+O(h^4),$ and ignoring higher order terms, the equation \eqref{IQ3} become
\begin{equation}
u_{n+1}=( 1-\epsilon_n^2h^2)\bigg(\frac{h(1+\epsilon_n^2h^2)(2+\epsilon_n^2h^2)f_n+2u_n}{2}\bigg).
\end{equation}
From the local truncation error at point  $t=t_n$, we get
\begin{equation}
\tau_n=\bigg(u_n\epsilon_n^2+\frac{u_n''}{2}\bigg)h+\bigg(\frac{u_n'''}{6}-\frac{u_n'\epsilon_n^2}{2}\bigg)h^2+O(h^3).
\end{equation}
For $O(h^2)$,
\begin{equation}
\epsilon_n^2=-\frac{u_n''}{2u_n}.
\end{equation}
Given $u_n'=f_n$  and therefore $u_n''=f_n'$. Using backward difference formula for $f_n'=\displaystyle \dfrac{f_n-f_{n-1}}{h}$, we obtain optimal shape parameter as
\begin{equation}
\epsilon_n^2=-\dfrac{(f_n-f_{n-1})}{2hu_n}.
\end{equation}
Thus, we have  adaptive modified IQ-RBF Euler method 
\begin{equation}
u_{n+1}=\displaystyle \bigg(1+\dfrac{(f_n-f_{n-1})h}{2u_n}\bigg)\bigg(\dfrac{h}{2}\bigg(1-\dfrac{(f_n-f_{n-1})h}{2u_n}\bigg)\bigg(2-\dfrac{(f_n-f_{n-1})h}{2u_n}\bigg)f_n+u_n\bigg).
\end{equation}
To avoid the denominator become zero, the condition C1 in \eqref{C1} is same for IQ-RBF whereas the condition C2 become
\begin{eqnarray}\label{C12}
\text{C2}: \epsilon_n^2=
\begin{cases}
\text{sgn} \bigg((f_n-f_{n-1}) 2hu_n \bigg)  L  & \text{if}\,\,\,|u_n| \leq h^p,\\
-\dfrac{f_n-f_{n-1}}{2hu_n} & \text{otherwise.}\,\,\,\, 
\end{cases}
\end{eqnarray}

\section{Stability and Convergence}
\label{sec: SC}
In this section, we provide the convergence and stability of proposed RBF Euler method. In \cite{Dah}, it is shown that the linear method is convergent if it is consistent and convergence. However, the proposed methods are not linear, thus we cannot
apply this technique but need to explicitly show the convergence and  stability of each derived method. 

\subsection{Adaptive IMQ-RBF Euler method}
\label{sec: SIMQRBF}
The adaptive IMQ-RBF Euler method \eqref{IMQE1-RBF} takes the form
\begin{equation*}
\displaystyle u_{n+1}=\frac{(1+\epsilon_n^2h^2)hf_n+u_n}{\sqrt{1+\epsilon_n^2h^2}}.
\end{equation*}
On simplifying, we get
\begin{equation*}
\displaystyle u_{n+1}=\sqrt{(1+\epsilon_n^2h^2)}hf_n+\frac{u_n}{\sqrt{1+\epsilon_n^2h^2}},
\end{equation*}
and the truncation error is
\begin{equation*}
\displaystyle\tau_n=\displaystyle \dfrac{u_{n+1}-\dfrac{(1+\epsilon_n^2h^2)hf_n+u_n}{\sqrt{1+\epsilon_n^2h^2}}}{h}.
\end{equation*}
The exact solution satisfies
\begin{equation*}
\displaystyle v_{n+1}=\frac{(1+\epsilon_n^2h^2)hf_n+v_n}{\sqrt{1+\epsilon_n^2h^2}}+h\tau_n.
\end{equation*}
On simplification, we get
\begin{equation*}
\displaystyle v_{n+1}=\sqrt{(1+\epsilon_n^2h^2)}hf_n+\frac{v_n}{\sqrt{1+\epsilon_n^2h^2}}+h\tau_n.
\end{equation*}
Subtracting $u_{n+1}$ from $v_{n+1}$, we obtain
\begin{eqnarray*}
\begin{aligned}
E_{n+1}= & \lvert v_{n+1}-u_{n+1}\rvert\\
= & \bigg\lvert \frac{1}{\sqrt{1+\epsilon_n^2h^2}}(v_n-u_n)+h\sqrt{1+\epsilon_n^2h^2}[f(t_n,v_n)-f(t_n,u_n)]+h\tau_n \bigg\rvert,\\
 \leq &  \frac{1}{\sqrt{1+\epsilon_n^2h^2}}\lvert(v_n-u_n)\rvert+h\sqrt{1+\epsilon_n^2h^2}\lvert f(t_n,v_n)-f(t_n,u_n)\rvert+h\lvert\tau_n\rvert,\\
\leq & \frac{1}{\sqrt{1+\epsilon_n^2h^2}} E_n+hL_f\sqrt{1+\epsilon_n^2h^2}\lvert v_n-u_n\rvert+h\lvert\tau_n\rvert,\\
= & \frac{1}{\sqrt{1+\epsilon_n^2h^2}}\displaystyle (1+h(1+\epsilon_n^2h^2)L_f)E_n+h\lvert \tau_n \rvert.
\end{aligned}
\end{eqnarray*}
Here, we assume that the function  $f(t, u)$ is  Lipschitz continuious in $u$,
\begin{equation*}
\lvert f(t_n,v_n)-f(t_n,u_n) \rvert \leq L_f \lvert v_n -u_n \rvert = L_f E_n.
\end{equation*}
Using Induction process,  we can show that
\begin{eqnarray*}
\begin{aligned}
E_n \leq & \prod_{j=0}^{n-1} \frac{1}{\sqrt{1+\epsilon_j^2h^2}}\displaystyle (1+h(1+\epsilon_j^2h^2)L_f)E_0 \\
& +h\sum_{j=1}^{n-1}\prod_{m=j}^{n-1}\frac{1}{\sqrt{1+\epsilon_m^2h^2}}\displaystyle (1+h(1+\epsilon_m^2h^2)L_f)^{n-j}\lvert \tau_{j-1}\rvert +h\lvert \tau_{n-1}\rvert.
\end{aligned}
\end{eqnarray*}
Consider $N$ number of time steps needed to reach time b, that is, $N=(b-a)/h$ and set
 \begin{equation*}
{\lvert \lvert \tau \rvert \rvert}_{\infty} = \underset{0 \leq n \leq N-1}{\mathrm{max}} \lvert \tau_n \rvert.
\end{equation*}
Assuming initial condition is exact, we get $E_0=0$. Therefore, as $h \to 0$,
\begin{equation*}
E_n \leq e^{L_fT} T {\lvert \lvert \tau \rvert \rvert }_{\infty}.
\end{equation*}
This is valid for all $n$ which satisfy the condition  $nh \leq T=b-a$. This proves that the adaptive IMQ-RBF Euler method is convergent.

\subsection{Adaptive IQ-RBF Euler method}
\label{sec: SIQRBF}
The adaptive IQ-RBF Euler method \eqref{IQ3} takes the form
\begin{equation*}
\displaystyle u_{n+1}=\frac{(1+\epsilon_n^2h^2)(2+\epsilon_n^2h^2)hf_n+2u_n}{2(1+\epsilon_n^2h^2)}.
\end{equation*}
On further simplification, we get
\begin{equation*}
\displaystyle u_{n+1}=\frac{(2+\epsilon_n^2h^2)hf_n}{2}+\frac{u_n}{(1+\epsilon_n^2h^2)},
\end{equation*}
and the truncation error is
\begin{equation*}
\displaystyle\tau_n=\displaystyle \dfrac{u_{n+1}-\dfrac{(1+\epsilon_n^2h^2)(2+\epsilon_n^2h^2)hf_n+2u_n}{2(1+\epsilon_n^2h^2)}}{h}.
\end{equation*}
The exact solution satisfies
\begin{equation*}
\displaystyle v_{n+1}=\frac{(1+\epsilon_n^2h^2)(2+\epsilon_n^2h^2)hf_n+2v_n}{2(1+\epsilon_n^2h^2)}+h\tau_n
\end{equation*}
On simplification, we get
\begin{equation*}
\displaystyle v_{n+1}=\frac{(2+\epsilon_n^2h^2)hf_n}{2}+\frac{v_n}{(1+\epsilon_n^2h^2)}+h\tau_n.
\end{equation*}
Subtracting $u_{n+1}$ from $v_{n+1}$, we obtain
\begin{eqnarray*}
\begin{aligned}
E_{n+1}= & \lvert v_{n+1}-u_{n+1}\rvert\\
= &  \bigg\lvert \frac{1}{(1+\epsilon_n^2h^2)}(v_n-u_n)+\frac{h(2+\epsilon_n^2h^2)}{2}[f(t_n,v_n)-f(t_n,u_n)]+h\tau_n\bigg\rvert, \\
 \leq & \frac{1}{(1+\epsilon_n^2h^2)}\lvert(v_n-u_n)\rvert+\frac{h(2+\epsilon_n^2h^2)}{2}\lvert f(t_n,v_n)-f(t_n,u_n)\rvert+h\lvert\tau_n\rvert, \\
\leq &  \frac{1}{(1+\epsilon_n^2h^2)} E_n+hL_f\frac{(2+\epsilon_n^2h^2)}{2}\lvert v_n-u_n\rvert+h\lvert\tau_n\rvert, \\
= & \bigg(\frac{1}{2(1+\epsilon_n^2h^2)}\displaystyle (2+h(2+\epsilon_n^2h^2)(1+\epsilon_n^2h^2)L_f)\bigg)E_n+h\lvert \tau_n \rvert.
\end{aligned}
\end{eqnarray*}
By induction principlie,
\begin{eqnarray*}
\begin{aligned}
E_n \leq & \prod_{j=0}^{n-1} \frac{1}{2(1+\epsilon_n^2h^2)}\displaystyle (2+h(1+\epsilon_n^2h^2)(2+\epsilon_n^2h^2)L_f)E_0\\
& +h\sum_{j=1}^{n-1}\prod_{m=j}^{n-1}\frac{1}{2(1+\epsilon_m^2h^2)}\displaystyle (2+h(1+\epsilon_m^2h^2)(2+\epsilon_m^2h^2)L_f)^{n-j}\lvert \tau_{j-1}\rvert +h\lvert \tau_{n-1}\rvert.
\end{aligned}
\end{eqnarray*}
Set ${\lvert \lvert \tau \rvert \rvert}_{\infty}=\underset{0 \leq n \leq N-1}{\mathrm{max}} \lvert \tau_n \rvert, $ $E_0=0$ and as $h \to 0$, we get 
\begin{equation*}
E_n \leq e^{L_fT} T {\lvert \lvert \tau \rvert \rvert }_{\infty},
\end{equation*}
 which is valid for all $n$ with  $nh \leq T.$ This proves that adaptive IQ-RBF Euler method is convergent.
\par
\begin{figure}[!h]
\begin{center}
\includegraphics[trim=45mm  10mm 180mm 10mm, clip=true, scale=0.4]{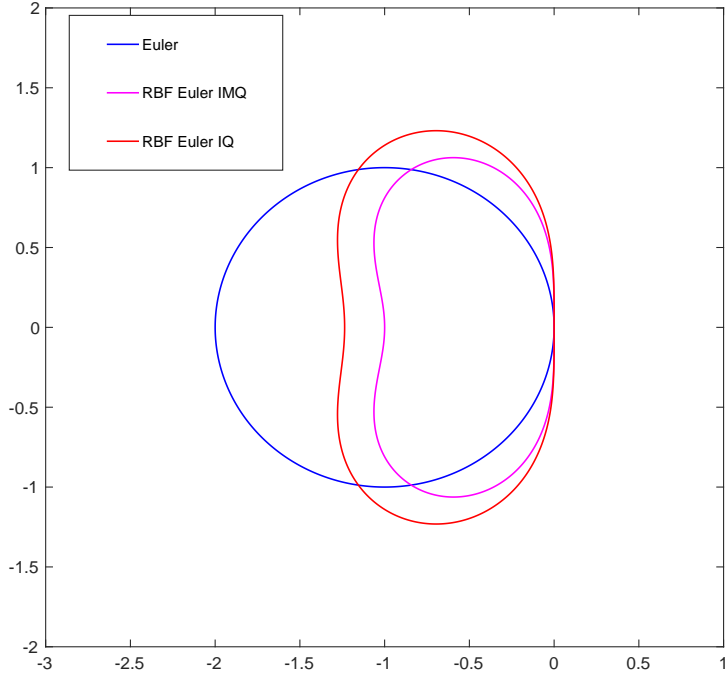}
\caption{Stability region for Euler, IQ- and IMQ-RBF Euler methods}\label{S1}
\end{center}
\end{figure}

The Figure \eqref{S1} shows the stability region for Euler, IQ-RBF Euler  and IMQ-RBF Euler methods. 
In the following,  we collect the Euler and adaptive Euler methods developed in \cite{GJ-1, GJ-2} and the methods proposed here with optimal shape parameters.
\begin{table}[ht!]
\label{Table:i}
\caption{Collection of   adaptive Euler methods}
\begin{center}
\scriptsize
\begin{tabular}{ l | l | l | p{2.5cm} }
\hline 
\textbf{Method} & \textbf{Numerical scheme} & \textbf{Order} & \textbf{Optimal $\epsilon_n^2$} \\
\hline
 \textit{Euler}  & $u_{n+1}=u_n+hf_n$  & $O(h)$ & \\ 
\textit{MQ-RBF Euler}  & $u_{n+1}=\bigg(1+\dfrac{\epsilon_n^2h^2}{2}\bigg)(u_n+hf_n)$ & $O(h^2)$ & $\dfrac{f_n-f_{n-1}}{hu_n}$\\
\textit{Gaussian-RBF Euler}  & $u_{n+1}=u_ne^{-\epsilon_n^2h^2}+hf_n$ & $O(h^2)$ & $-\dfrac{f_n-f_{n-1}}{2hu_n}$\\
 \textit{IMQ-RBF Euler}  & $u_{n+1}=\dfrac{(1+\epsilon_n^2h^2)hf_n+u_n}{\sqrt{1+\epsilon_n^2h^2}}$ & $O(h^2) $& $-\dfrac{f_n-f_{n-1}}{hu_n}$\\
 \textit{IQ-RBF Euler}  & $u_{n+1}=\dfrac{h(1+\epsilon_n^2h^2)(2+\epsilon_n^2h^2)f_n+2u_n}{2(1+\epsilon_n^2h^2)}$ & $O(h^2) $& $-\dfrac{f_n-f_{n-1}}{2hu_n}$\\
\hline
\end{tabular}
\end{center}
\end{table}

\section{Numerical experiments}
\label{sec:nr}
In this section, we provide  four examples to verify the  numerical results introduced in the previous sections.
\subsection*{Example 1}
\begin{figure}[!h]
\begin{center}
\includegraphics[trim=0cm  10mm 0.2cm 0cm, clip=true, , scale=0.32]{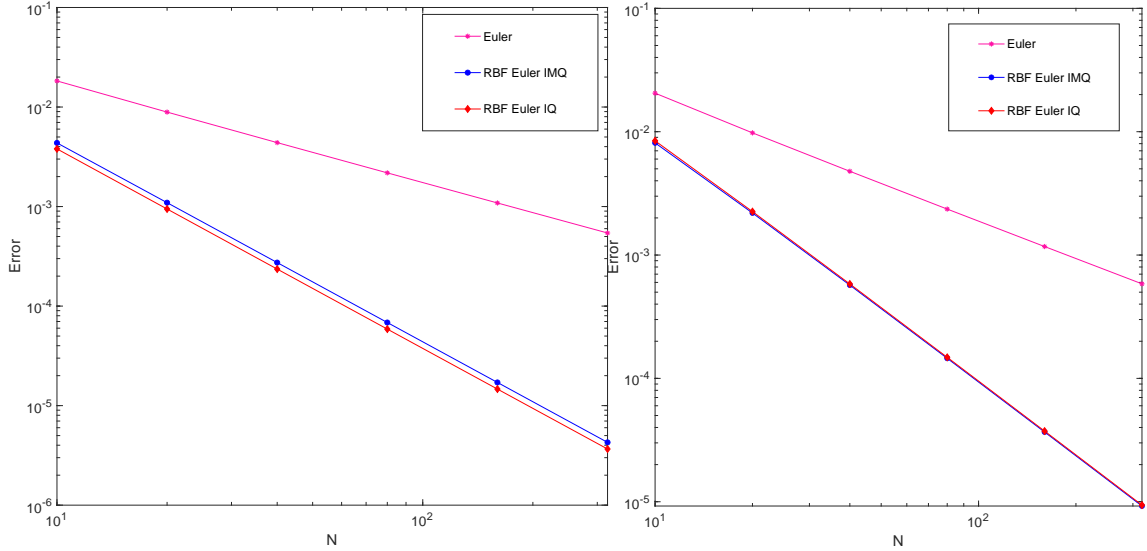}
\caption{The global errors versus $N$ in logarithmic scale for example-\ref{ex1} (left) and \ref{ex2} (right)}
\label{fig:2}
\end{center}
\end{figure}

\normalfont
We consider the following Initial value problem
\begin{equation}\label{ex1}
\frac{du}{dt}=-u^2,0<t \leq 1, u(0)=1.
\end{equation}
Table \ref{Table:1} shows the global errors versus various N and local convergence orders.  From this table \ref{Table:1}, it is observed that Euler method is of first order accurate whereas adaptive Euler methods achieve second order of accuracy. Figure \ref{fig:2} (left) shows the global errors versus $N$ at $t=1$ by the Euler’s method,  the Euler method with Inverse-Quadratic and Inverse Multi-Quadratic  RBFs. It concludes that from the figure \ref{fig:2} (left), the Euler’s method yields the first order accuracy while the RBF Euler’s methods yield the second order accuracy.
\begin{table}[ht!]
\label{Table:1}
\caption{ Global errors at $t=1$ and order of convergence for example (\ref{ex1})}
\begin{center}
\scriptsize
\begin{tabular}{|c c c c|}
\hline
Method & N & Global Error & Order\\ [0.5 ex]
\hline
Euler &10 & 0.018287121529848 & -  	\\
&20 & 0.008895076334408 & 1.0397\\
&40 & 0.004388827380214 & 1.0192\\
&80 & 0.002180125588386 & 1.0094\\
&160 & 0.001086537438631 & 1.0047\\
&320 & 0.000542393094490 & 1.0023\\[1 ex]
\hline
RBF Euler IMQ& 10 &	      0.004359450155230   &-  	\\[0.5 ex]
&20&	 0.001093900148224 &  	1.9947  \\ 
&40	& 0.000273574228572 &  	1.9995 \\  
&80	&0.000068383354031  & 	2.0002 \\  
&160&	0.000017093255390 &  	2.0002 \\  
&320&	0.000004272912760&	2.0001\\[1 ex]
\hline
RBF Euler IQ& 10&	0.003796432501710     &-  	\\[0.5 ex]
&20&	0.000944949927189  & 	2.0063     \\ 
&40&	0.000235395279817  & 	2.0052    \\  
&80&	0.000058726355740  & 	2.0030    \\  
&160&	0.000014665313947  & 	2.0016    \\  
&320&	0.000003664237210	&2.0008\\[1 ex]
\hline
\end{tabular}
\end{center}
\end{table}
\begin{table}[ht!]
\label{Table:2}
\caption{ Global errors at $t=2$ and order of convergence for example (\ref{ex2})}
\begin{center}
\scriptsize
\begin{tabular}{|c c c c|}
\hline
Method & N & Global Error & Order\\ [0.5 ex]
\hline
Euler &10&	0.020530545634685 &-\\ 	
&20&	0.009786516161729&   	1.0689  \\
&40&	0.004775483993213 &  	1.0351   \\
&80&	0.002358560616812&   	1.0177   \\
&160	&0.001172019387861  & 	1.0089   \\
&320	&0.000584198876871&	1.0045\\[1 ex]
\hline
RBF Euler IMQ&10&	0.008113825603093 &-\\[0.5 ex]	
&20&	0.002192315187090 &  	1.8879    \\
&40&	0.000570268629560&   	1.9427   \\
&80&	0.000145465457719 &  	1.9710 \\
&160&	0.000036737092202&   	1.9854   \\
&320&	0.000009231163763&	1.9927\\[1 ex]
\hline
RBF Euler IQ&10	&0.008429894186101   &-\\[0.5 ex]		
&20&	0.002251112343783 &  	1.9049       \\
& 40&	0.000582515300086 &  	1.9503   \\
&80&	0.000148226787150 &  	1.9745   \\
&160&	0.000037390349041  & 	1.9871     \\
&320&	0.000009389875731	&1.9935\\[1 ex]
\hline
\end{tabular}
\end{center}
\end{table}
\subsection*{Example 2}
\normalfont
Let us consider the second problem 
\begin{equation}\label{ex2}
u'=\frac{2t^2-u}{t^2u-t},1<t \leq 2, u(1)=2,
\end{equation}
which is in  non separable problem.   The exact solution to the differential equation is $u(t)=\frac{1}{t}+\sqrt{\frac{1}{t^2}+4t-4}.$  The  figure \ref{fig:2}(right) shows the global errors versus $N$ at $t=2$ by the  Euler’s method,  Inverse-Quadratic RBF Euler and Inverse Multi-Quadratic RBF Euler method.  It is shown in the figure, the Euler’s method yields the first order accuracy while the RBF Euler’s methods performs well and yield the second order accuracy. Table \ref{Table:2} shows the global errors versus various N and local convergence orders.
 \begin{figure}[ht!]
\begin{center}
\includegraphics[trim=8mm  1mm 0.2cm 0cm, clip=true, , scale=0.32]{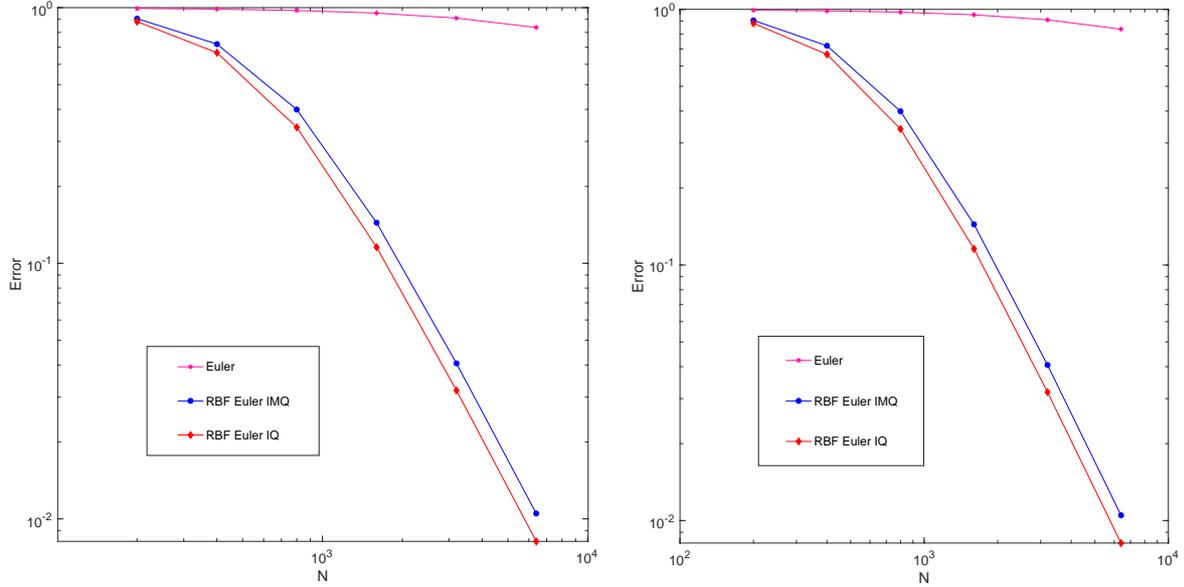}
\caption{The global errors versus $N$ in logarithmic scale for example (\ref{ex3})(left) and (\ref{ex4}) with no condition (right)}
\label{fig:4}
\end{center}
\end{figure}
\begin{table}[ht!]
\label{Table:3}
\caption{ Global errors at $t=0$ and order of convergence for example (\ref{ex3})}
\begin{center}
\scriptsize
\begin{tabular}{|c c c c|}
\hline
Method & N & Global Error & Order\\ [0.5 ex]
\hline
Euler&200	&0.992928300529281    &-\\	
&400&	0.986794866203422  & 	0.0089     \\
&800&	0.974934408048963&   	0.0175      \\
&1600&	0.952262436431508&   	0.0339     \\
&3200&	0.910174769390209&   	0.0652    \\
&6400&	0.836481000589044&	0.1218\\[1 ex]
\hline
RBF Euler IMQ&200&	0.905715491232480&    	-\\[0.5 ex]	
&400&	0.720151335458197&   	0.3308      \\
&800&	0.399219472430586&   	0.8511     \\
&1600&	0.144107798700244&   	1.4700    \\
&3200&	0.040598570535400&   	1.8276      \\
&6400&	0.010491065260367&	1.9523\\[1 ex]
\hline
RBF Euler IQ&200&	0.881748559214363&   	-\\[0.5 ex]	
&400&	0.666390343964605&   	0.4040    \\
&800&	0.340232327363573 &  	0.9698     \\
&1600&	0.115568657862118 &  	1.5578    \\
&3200&	0.031803564115116 &  	1.8615         \\
&6400&	0.008164385861870&	1.9618\\[1 ex]
\hline
\end{tabular}
\end{center}
\end{table}
\begin{figure}[ht!]
\begin{center}
\includegraphics[trim=2.5cm  0.5cm 0.2cm 0.2cm, clip=true, , scale=0.35]{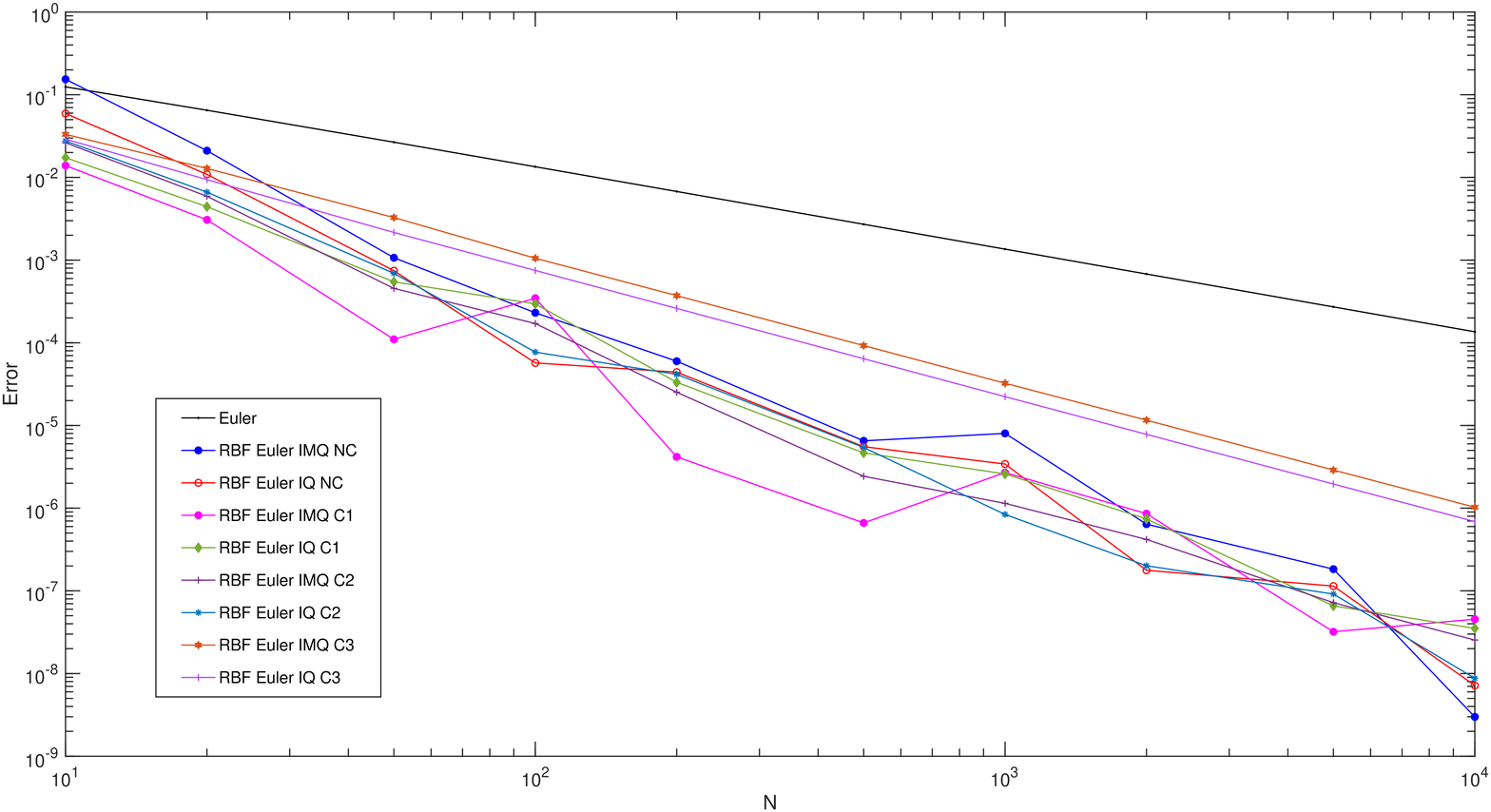}
\caption{The global errors versus $N$ in logarithmic scale for example (\ref{ex4}) with different conditions}
\label{fig:5}
\end{center}
\end{figure}
\subsection*{Example 3}
Let us consider the problem
\begin{equation}\label{ex3}
u'=-4t^3u^2,-10<t \leq 0,
u(-10)=\dfrac{1}{10001},
\end{equation}
which is again in non separable problem where the solution of this problem rapidly changes. The exact solution to the differential equation is $u(t)=\frac{1}{t^4+1}.$ We choose this example to show how the proposed RBF methods behave for a stiff problem. The left part of Figure \ref{fig:4} shows that the global errors versus $N$ at $t=0$ by the Euler’s method, the RBF Euler method with Inverse Multi-Quadratic and  Inverse-Quadratic functions.   From the Table \ref{Table:3}, it is observed that Euler method is not even achieving the first order accuracy whereas adaptive Euler methods are close to second order convergence as mesh size become small. 
\subsection*{Example 4}
Consider the following differential equation
\begin{equation}\label{ex4}
\frac{du}{dt}=u+2,0<t \leq 1,
\end{equation}
with initial condition
\begin{equation}
u(0)=-1.
\end{equation}
\begin{table}[ht!]
\label{Table:4}
\caption{ Global errors at $t = 1$ and orders of convergence for example (\ref{ex4}) with no condition}
\begin{center}
\scriptsize
\begin{tabular}{|c c c c|}
\hline
Method & N & Global Error & Order\\ [0.5 ex]
\hline
Euler &10 &   0.124539368359045 & -  	\\
&20 & 0.064984123314625    &     0.9384   \\
&50 & 0.026693799385440    & 1.2836   \\
&100 & 0.013467999037519    & 0.9870   \\
&200 & 0.006764705529671    & 0.9934   \\
&500 & 0.002713307807321 & 1.3180\\
&1000 &  0.001357896223155&   0.9987     \\
 &2000 & 0.000679259135906 & 0.9993     \\
&5000 & 0.000271778357187       & 1.3215      \\
&10000 & 0.000135901633822& 0.9999   \\[1 ex]
\hline
RBF Euler IMQ& 10 & 0.153968657361701 &-	\\
&20&	 0.021060740746961    &  	2.8700     \\ 
&50	& 0.001068383514315  &  	4.3011    \\  
&100&	0.000230835080698     & 	2.2105    \\  
&200&	0.000059859695074    &  	1.9472    \\  
&500&	0.000006519915393&	3.1987\\ 
&1000 &  0.000008021749547   &   -0.2991        \\
 &2000 & 0.000000640575256    & 3.6465        \\
&5000 &  0.000000182984119         & 1.8076         \\
&10000 & 0.000000002989985&  5.9354 \\[1 ex]
\hline
RBF Euler IQ& 10&	0.059308599951798       &-  	\\
&20&	 0.010880410041919    & 	2.4465        \\ 
&50&	0.000744746765859     & 	3.8688       \\  
&100&	0.000057093598870     & 	3.7053       \\  
&200&	0.000044183937295     & 	0.3698       \\  
&500&	0.000005545405298	&2.9942\\
&1000 &  0.000003417146218 &0.6985         \\
 &2000 &0.000000177444964       & 4.2673    \\
&5000 &  0.000000114231586  & 0.6354    \\
&10000 & 0.000000007166101&  3.9946 \\[1 ex]
\hline
\end{tabular}
\end{center}
\end{table}
\begin{table}[ht!]
\label{Table:5}
\caption{ Global errors at $t = 1$ and orders of convergence for example (\ref{ex4}) with condition p=1,L=0}
\begin{center}
\scriptsize
\begin{tabular}{|c c c c|}
\hline
Method & N & Global Error & Order\\ [0.5 ex]
\hline
RBF Euler IMQ& 10&	0.013926573217714   &-	\\[0.5 ex]
&20&0.003073329603225  & 	2.1800   \\ 
&50&	0.000110087436249   & 	4.8031     \\  
&100&	0.000346251852890    & 	-1.6532  \\  
&200&	0.000004173836753     &	6.3743    \\  
&500&	0.000000663116536	&2.6540\\
&1000 &  0.000002696484348    &2.0237   \\
 &2000 &0.000000857921772    & 1.6522       \\
&5000 &  0.000000031988806    & 4.7452      \\
&10000 &0.000000045408762&  -0.5054 \\[1 ex]
\hline
RBF Euler IQ& 10&	0.017344120771895 &-  \\[0.5 ex]
&20&	 0.004428650824286      & 1.9695      \\ 
&50&	0.000546495750087        & 	3.0186     \\  
&100&	   0.000295540864736        & 0.8869   \\  
&200&	 0.000033266543641      & 	3.1512    \\  
&500&	0.000004662448778	&2.8349\\
&1000 &    0.000002593720385    & 0.8461    \\
 &2000 &0.000000740408153  &    1.8086    \\
&5000 & 0.000000065844749     & 3.4912       \\
&10000 &0.000000035164452 &   0.9050\\[1 ex]
\hline
\end{tabular}
\end{center}
\end{table}
\begin{table}[ht!]
\label{Table:6}
\caption{ Global errors at $t=1$ and orders of convergence for example (\ref{ex4}) with condition p=1, $L=1/h$}
\begin{center}
\scriptsize
\begin{tabular}{|c c c c|}
\hline
Method & N & Global Error & Order\\ [0.5 ex]
\hline
RBF Euler IMQ& 10&  0.026264826889257  &- \\[0.5 ex]
&20&	 0.005915299559272  & 	2.1506     \\ 
&50&	0.000455177536653   & 	3.6999    \\  
&100&	0.000170920081854   & 	1.4131    \\  
&200&	0.000025243414411   & 	2.7593     \\  
&500&	0.000002429153891	&3.3774\\
&1000 &  0.000001142798939  &   1.0879     \\
 &2000 &0.000000418791323    & 1.4483      \\
&5000 &  0.000000071959935    & 2.5410      \\
&10000 &0.000000025409851 &  1.5018 \\[1 ex]
\hline
RBF Euler IQ& 10&	0.017344120771895 &-\\[0.5 ex]
&20&	    0.004428650824286   & 	1.9695     \\ 
&50&	0.000546495750087   & 3.0186      \\  
&100&	0.000295540864736    & 	0.8869    \\  
&200&	0.000033266543641  & 	3.1512    \\  
&500&	0.000004662448778	&2.8349\\
&1000 &    0.000002593720385    &   0.8461 \\
 &2000 & 0.000000740408153   &1.8086      \\
&5000 &  0.000000065844749     &   3.4912  \\
&10000 &0.000000035164452&  0.9050 \\[1 ex]
\hline
\end{tabular}
\end{center}
\end{table}
\begin{table}[ht!]
\label{Table:7}
\caption{ Global errors at $t=1$ and orders of convergence for example (\ref{ex4}) with condition p=1, $L=1/h^{0.5}$}
\begin{center}
\scriptsize
\begin{tabular}{|c c c c|}
\hline
Method & N & Global Error & Order\\ [0.5 ex]
\hline
RBF Euler IMQ& 10&	0.033178146726859 &-  	\\[0.5 ex]
&20&	   0.012849172714074   &  1.3686    \\ 
&50&	0.003275119747404   & 	1.9721     \\  
&100&	0.001053900009690  & 	1.6358     \\  
&200&	0.000372028816605  & 	1.5023     \\  
&500&	0.000092572462070	&2.0068\\
&1000 &  0.000032403967052    &   1.5144  \\
 &2000 &   0.000011602143029 & 1.4818      \\
&5000 &  0.000002890485356 &2.0050   \\
&10000 &  0.000001025943622   &  1.4944 \\[1 ex]
\hline
RBF Euler IQ& 10&  0.029061309514811 &-\\[0.5 ex]
&20&	 0.009419481977583  & 	 1.6254        \\ 
&50&	 0.002159530710872  &  2.1249       \\  
&100&	0.000749695742409   & 1.5263       \\  
&200&	 0.000260534051296  & 	1.5248     \\  
&500&	0.000064144990861	&  2.0220    \\
&1000 &     0.000022316165521   &  1.5232  \\
 &2000 &0.000007810878018   &1.5145    \\
&5000 &  0.000001957655265  &1.9963    \\
&10000 & 0.000000688491352&  1.5076 \\[1 ex]
\hline
\end{tabular}
\end{center}
\end{table}

The exact solution to the differential equation is $u(t)=e^t-2$ and $u(t)=0$ at $t=\ln 2 \approx 0.6931.$   Table \ref{Table:4} and  Figure \ref{fig:4}(right) shows the global errors versus $N$ at $t=1$ by the Euler’s method,  the  proposed RBF Euler method with Inverse Multi-Quadratic and Inverse-Quadratic functions. Tables \ref{Table:4}, \ref{Table:5}, \ref{Table:6} and \ref{Table:7}  shows the errors and orders for the cases; 
\begin{itemize}
\item NC: no condition,
\item  C1: with $p=1$, $L=0$ in \eqref{C2}, \eqref{C12},
\item  C2:  with $p=1$, $L=1/h$ in \eqref{C2}, \eqref{C12}, and
\item C3: with $p=1$, $L=1/h^{0.5}$ in \eqref{C2}, \eqref{C12}. 
\end{itemize}
 Figure \ref{fig:5} shows the global errors versus $N$ at $t=1$ for the original Euler’s method with the  proposed Inverse Multi-Quadratic and Inverse-Quadratic RBF Euler method with the NC, C1, C2 and C3 conditions.  It is observed from the figure that all the RBF methods are stable, they all yield the overall second order convergence. Particularly, if  $p=1,$ and $L=1/h^{0.5}$, both the schemes achieve second order convergence rate and the same is shown in  figure \ref{fig:5}.
\section{Conclusion}
\label{sec:con}
In this paper, we proposed two enhanced Euler methods to solve initial value problems.  These methods seek the approximation using inverse-quadratic and inverse multi-quadratic  radial basis functions interpolation. The RBF functions have a free shape parameter and by optimising this, the modified Euler methods yields more accurate results in comparison to the classical Euler  ODE solver. The consistency, convergence analysis and stability regions are studied.  This is the preliminary study of the shape parameters for the  inverse-quadratic and inverse multi-quadratic radial basis functions and extension to the multi-quadratic and Gaussian-RBFs \cite{GJ-1, GJ-2}. Further research is going on that, how to choose the best performing shape parameter and how the modified methods do work for more complicated ODEs and PDEs.

\end{document}